\documentclass[12pt,oneside]{amsart}
\usepackage[T1]{fontenc}

\numberwithin{equation}{section}
\newtheorem{thm}[equation]{Theorem}

\newtheorem{lem}[equation]{Lemma}

\newenvironment{pf}{\proof[\proofname]}{\endproof}
\newenvironment{pf*}[1]{\proof[#1]}{\endproof}

\theoremstyle{definition}

\theoremstyle{remark}

\newtheorem*{ack}{Acknowledgement}

\usepackage{amsmath,amsthm}
\usepackage{amsfonts}

\makeatother

\begin{document}
\baselineskip=18truept


\def\C {{\mathbb C}}
\def\Cn {{\mathbb C}^n}
\def\R {{\mathbb R}}
\def\Rn {{\mathbb R}^n}
\def\Z {{\mathbb Z}}
\def\N {{\mathbb N}}
\def\cal#1{{\mathcal #1}}
\def\bb#1{{\mathbb #1}}

\def\dbar {\bar \partial }
\def\dir {{\mathcal D}}
\def\lev#1{{\mathcal L}(#1)}
\def\lap {\Delta }
\def\ol {{\mathcal O}}
\def\E {{\mathcal E}}
\def\J {{\mathcal J}}
\def\U {{\mathcal U}}
\def\V {{\mathcal V}}
\def\z {\zeta }
\def\Harm {\text {Harm}\, }
\def\grad {\nabla }
\def\dexh {\{ M_k \} _{k=0}^{\infty } }
\def\sing#1{#1_{\text {sing}}}
\def\reg#1{#1_{\text {reg}}}

\def\setof#1#2{\{ \, #1 \mid #2 \, \} }
\def\Subset{\subset\subset }

\def\holecl {M\setminus \overline M_0}
\def\hole {M\setminus M_0}

\def\nd{\frac {\partial }{\partial\nu } }
\def\ndof#1{\frac {\partial#1}{\partial\nu } }

\def\pdof#1#2{\frac {\partial#1}{\partial#2}}

\def\cinf{C^{\infty }}

\def\diam{\text {diam} \, }

\def\real{\text {Re}\, }

\def\imag{\text {Im}\, }

\def\supp{\text {supp}\, }

\def\Vol{\text {\rm vol} \, }



\def\anal{analytic }
\def\analns{analytic}

\def\bdd{bounded }
\def\bddns{bounded}

\def\cpt{compact }
\def\cptns{compact}

\def\cpx{complex }
\def\cpxns{complex}

\def\cont{continuous }
\def\contns{continuous}

\def\dime{dimension }
\def\dimens{dimension }

\def\exh{exhaustion }
\def\exhns{exhaustion}

\def\fn{function }
\def\fnns{function}

\def\fns{functions }
\def\fnsns{functions}

\def\holo{holomorphic }
\def\holons{holomorphic}

\def\mero{meromorphic }
\def\merons{meromorphic}

\def\holoconvex{holomorphically convex }
\def\holoconvexns{holomorphically convex}

\def\ircomp{irreducible component }
\def\concomp{connected component }
\def\ircompns{irreducible component}
\def\concompns{connected component}
\def\ircomps{irreducible components }
\def\concomps{connected components }
\def\ircompsns{irreducible components}
\def\concompsns{connected components}

\def\irred{irreducible }
\def\irredns{irreducible}

\def\con{connected }
\def\conns{connected}

\def\comp{component }
\def\compns{component}
\def\comps{components }
\def\compsns{components}

\def\mfld{manifold }
\def\mfldns{manifold}
\def\mflds{manifolds }
\def\mfldsns{manifolds}

\def\nbd{neighborhood }
\def\nbds{neighborhoods }
\def\nbdns{neighborhood}
\def\nbdsns{neighborhoods}

\def\harm{harmonic }
\def\harmns{harmonic}
\def\plh{pluriharmonic }
\def\plhns{pluriharmonic}
\def\plsh{plurisubharmonic }
\def\plshns{plurisubharmonic}

\def\qplsh#1{$#1$-plurisubharmonic}
\def\hplsh{$(n-1)$-plurisubharmonic }
\def\hplshns{$(n-1)$-plurisubharmonic}

\def\para{parabolic }
\def\parans{parabolic}

\def\rel{relatively }
\def\relns{relatively}

\def\str{strictly }
\def\strns{strictly}

\def\strg{strongly }
\def\strgns{strongly}

\def\cvx{convex }
\def\cvxns{convex}

\def\wrt{with respect to }
\def\wrtns{with respect to}

\def\st {such that }
\def\stns {such that}

\def\hm {harmonic measure }
\def\hmns {harmonic measure}

\def\hmib {harmonic measure of the ideal boundary of }
\def\hmibns {harmonic measure of the ideal boundary of}


\def\atil{\tilde a}
\def\btil{\tilde b}
\def\ctil{\tilde c}
\def\dtil{\tilde d}
\def\etil{\tilde e}
\def\ftil{\tilde f}
\def\gtil{\tilde g}
\def\htil{\tilde h}
\def\itil{\tilde i}
\def\jtil{\tilde j}
\def\ktil{\tilde k}
\def\ltil{\tilde l}
\def\mtil{\tilde m}
\def\ntil{\tilde n}
\def\otil{\tilde o}
\def\ptil{\tilde p}
\def\qtil{\tilde q}
\def\rtil{\tilde r}
\def\stil{\tilde s}
\def\ttil{\tilde t}
\def\util{\tilde u}
\def\vtil{\tilde v}
\def\wtil{\tilde w}
\def\xtil{\tilde x}
\def\ytil{\tilde y}
\def\ztil{\tilde z}

\def\Atil{\tilde A}
\def\Btil{\widetilde B}
\def\Ctil{\widetilde C}
\def\Dtil{\widetilde D}
\def\Etil{\widetilde E}
\def\Ftil{\widetilde F}
\def\Gtil{\widetilde G}
\def\Htil{\widetilde H}
\def\Itil{\tilde I}
\def\Jtil{\widetilde J}
\def\Ktil{\widetilde K}
\def\Ltil{\widetilde L}
\def\Mtil{\widetilde M}
\def\Ntil{\widetilde N}
\def\Otil{\widetilde O}
\def\Ptil{\widetilde P}
\def\Qtil{\widetilde Q}
\def\Rtil{\widetilde R}
\def\Stil{\widetilde S}
\def\Ttil{\widetilde T}
\def\Util{\widetilde U}
\def\Vtil{\widetilde V}
\def\Wtil{\widetilde W}
\def\Xtil{\widetilde X}
\def\Ytil{\widetilde Y}
\def\Ztil{\widetilde Z}

\def\alphatil {\tilde \alpha  }
\def\betatil {\tilde \beta  }
\def\gammatil {\tilde \gamma  }
\def\deltatil {\tilde \delta }
\def\epsilontil {\tilde \epsilon  }
\def\varepsilontil {\tilde \varepsilon  }
\def\zetatil {\tilde \zeta  }
\def\etatil {\tilde \eta  }
\def\thetatil {\tilde \theta  }
\def\varthetatil {\tilde \vartheta  }
\def\iotatil {\tilde \iota  }
\def\kappatil {\tilde \kappa  }
\def\lambdatil {\tilde \lambda  }
\def\mutil {\tilde \mu  }
\def\nutil {\tilde \nu  }
\def\xitil {\tilde \xi  }
\def\pitil {\tilde \pi  }
\def\varpitil {\tilde \varpi  }
\def\rhotil {\tilde \rho  }
\def\varrhotil {\tilde \varrho  }
\def\sigmatil {\tilde \sigma  }
\def\varsigmatil {\tilde \varsigma  }
\def\tautil {\tilde \tau  }
\def\upsilontil {\tilde \upsilon  }
\def\phitil {\tilde \phi  }
\def\varphitil {\tilde \varphi  }
\def\chitil {\tilde \chi  }
\def\psitil {\tilde \psi  }
\def\omegatil {\tilde \omega  }

\def\Gammatil {\widetilde \Gamma  }
\def\Deltatil {\tilde \Delta }
\def\Thetatil {\widetilde \Theta  }
\def\Lambdatil {\tilde \Lambda  }
\def\Xitil {\widetilde \Xi  }
\def\Pitil {\widetilde \Pi  }
\def\Sigmatil {\widetilde \Sigma  }
\def\Upsilontil {\widetilde \Upsilon  }
\def\Phitil {\tilde \Phi  }
\def\Psitil {\widetilde \Psi  }
\def\Omegatil {\widetilde \Omega  }

\def\varGammatil {\widetilde \varGamma  }
\def\varDeltatil {\tilde \varDelta  }
\def\varThetatil {\widetilde \varTheta  }
\def\varLambdatil {\tilde \varLambda  }
\def\varXitil {\widetilde \varXi  }
\def\varPitil {\widetilde \varPi  }
\def\varSigmatil {\widetilde \varSigma  }
\def\varUpsilontil {\widetilde \varUpsilon  }
\def\varPhitil {\tilde \varPhi  }
\def\varPsitil {\widetilde \varPsi  }
\def\varOmegatil {\widetilde \varOmega  }

\def\boldGammatil {\widetilde \boldGamma  }
\def\boldDeltatil {\tilde \boldDelta  }
\def\boldThetatil {\widetilde \boldTheta  }
\def\boldLambdatil {\tilde \boldLambda  }
\def\boldXitil {\widetilde \boldXi  }
\def\boldPitil {\widetilde \boldPi  }
\def\boldSigmatil {\widetilde \boldSigma  }
\def\boldUpsilontil {\widetilde \boldUpsilon  }
\def\boldPhitil {\tilde \boldPhi  }
\def\boldPsitil {\widetilde \boldPsi  }
\def\boldOmegatil {\widetilde \boldOmega  }


\def\ahat{\hat a}
\def\bhat{\hat b}
\def\chat{\hat c}
\def\dhat{\hat d}
\def\ehat{\hat e}
\def\fhat{\hat f}
\def\ghat{\hat g}
\def\hhat{\hat h}
\def\ihat{\hat i}
\def\jhat{\hat j}
\def\khat{\hat k}
\def\lhat{\hat l}
\def\mhat{\hat m}
\def\nhat{\hat n}
\def\ohat{\hat o}
\def\phat{\hat p}
\def\qhat{\hat q}
\def\rhat{\hat r}
\def\shat{\hat s}
\def\that{\hat t}
\def\uhat{\hat u}
\def\vhat{\hat v}
\def\what{\hat w}
\def\xhat{\hat x}
\def\yhat{\hat y}
\def\zhat{\hat z}

\def\Ahat{\hat A}
\def\Bhat{\widehat B}
\def\Chat{\widehat C}
\def\Dhat{\widehat D}
\def\Ehat{\widehat E}
\def\Fhat{\widehat F}
\def\Ghat{\widehat G}
\def\Hhat{\widehat H}
\def\Ihat{\hat I}
\def\Jhat{\widehat J}
\def\Khat{\widehat K}
\def\Lhat{\widehat L}
\def\Mhat{\widehat M}
\def\Nhat{\widehat N}
\def\Ohat{\widehat O}
\def\Phat{\widehat P}
\def\Qhat{\widehat Q}
\def\Rhat{\widehat R}
\def\Shat{\widehat S}
\def\That{\widehat T}
\def\Uhat{\widehat U}
\def\Vhat{\widehat V}
\def\What{\widehat W}
\def\Xhat{\widehat X}
\def\Yhat{\widehat Y}
\def\Zhat{\widehat Z}

\def\alphahat {\hat \alpha  }
\def\betahat {\hat \beta  }
\def\gammahat {\hat \gamma  }
\def\deltahat {\hat \delta }
\def\epsilonhat {\hat \epsilon  }
\def\varepsilonhat {\hat \varepsilon  }
\def\zetahat {\hat \zeta  }
\def\etahat {\hat \eta  }
\def\thetahat {\hat \theta  }
\def\varthetahat {\hat \vartheta  }
\def\iotahat {\hat \iota  }
\def\kappahat {\hat \kappa  }
\def\lambdahat {\hat \lambda  }
\def\muhat {\hat \mu  }
\def\nuhat {\hat \nu  }
\def\xihat {\hat \xi  }
\def\pihat {\hat \pi  }
\def\varpihat {\hat \varpi  }
\def\rhohat {\hat \rho  }
\def\varrhohat {\hat \varrho  }
\def\sigmahat {\hat \sigma  }
\def\varsigmahat {\hat \varsigma  }
\def\tauhat {\hat \tau  }
\def\upsilonhat {\hat \upsilon  }
\def\phihat {\hat \phi  }
\def\varphihat {\hat \varphi  }
\def\vphihat {\hat \varphi  }
\def\chihat {\hat \chi  }
\def\psihat {\hat \psi  }
\def\omegahat {\hat \omega  }

\def\Gammahat {\widehat \Gamma  }
\def\Deltahat {\hat \Delta }
\def\Thetahat {\widehat \Theta  }
\def\Lambdahat {\hat \Lambda  }
\def\Xihat {\widehat \Xi  }
\def\Pihat {\widehat \Pi  }
\def\Sigmahat {\widehat \Sigma  }
\def\Upsilonhat {\widehat \Upsilon  }
\def\Phihat {\hat \Phi  }
\def\Psihat {\widehat \Psi  }
\def\Omegahat {\widehat \Omega  }

\def\varGammahat {\widehat \varGamma  }
\def\varDeltahat {\hat \varDelta  }
\def\varThetahat {\widehat \varTheta  }
\def\varLambdahat {\hat \varLambda  }
\def\varXihat {\widehat \varXi  }
\def\varPihat {\widehat \varPi  }
\def\varSigmahat {\widehat \varSigma  }
\def\varUpsilonhat {\widehat \varUpsilon  }
\def\varPhihat {\hat \varPhi  }
\def\varPsihat {\widehat \varPsi  }
\def\varOmegahat {\widehat \varOmega  }

\def\boldGammahat {\widehat \boldGamma  }
\def\boldDeltahat {\hat \boldDelta  }
\def\boldThetahat {\widehat \boldTheta  }
\def\boldLambdahat {\hat \boldLambda  }
\def\boldXihat {\widehat \boldXi  }
\def\boldPihat {\widehat \boldPi  }
\def\boldSigmahat {\widehat \boldSigma  }
\def\boldUpsilonhat {\widehat \boldUpsilon  }
\def\boldPhihat {\hat \boldPhi  }
\def\boldPsihat {\widehat \boldPsi  }
\def\boldOmegahat {\widehat \boldOmega  }

\def\seq#1#2{\{#1_{#2}\} }

\def\aseq{ \{ a_{\nu } \} }
\def\bseq{ \{ b_{\nu } \} }
\def\cseq{ \{ c_{\nu } \} }
\def\dseq{ \{ d_{\nu } \} }
\def\eseq{ \{ e_{\nu } \} }
\def\fseq{ \{ f_{\nu } \} }
\def\gseq{ \{ g_{\nu } \} }
\def\hseq{ \{ h_{\nu } \} }
\def\iseq{ \{ i_{\nu } \} }
\def\jseq{ \{ j_{\nu } \} }
\def\kseq{ \{ k_{\nu } \} }
\def\lseq{ \{ l_{\nu } \} }
\def\mseq{ \{ m_{\nu } \} }
\def\nseq{ \{ n_{\nu } \} }
\def\oseq{ \{ o_{\nu } \} }
\def\pseq{ \{ p_{\nu } \} }
\def\qseq{ \{ q_{\nu } \} }
\def\rseq{ \{ r_{\nu } \} }
\def\sseq{ \{ s_{\nu } \} }
\def\tseq{ \{ t_{\nu } \} }
\def\useq{ \{ u_{\nu } \} }
\def\vseq{ \{ v_{\nu } \} }
\def\wseq{ \{ w_{\nu } \} }
\def\xseq{ \{ x_{\nu } \} }
\def\yseq{ \{ y_{\nu } \} }
\def\zseq{ \{ z_{\nu } \} }

\def\Aseq{ \{ A_{\nu } \} }
\def\Bseq{ \{ B_{\nu } \} }
\def\Cseq{ \{ C_{\nu } \} }
\def\Dseq{ \{ D_{\nu } \} }
\def\Eseq{ \{ E_{\nu } \} }
\def\Fseq{ \{ F_{\nu } \} }
\def\Gseq{ \{ G_{\nu } \} }
\def\Hseq{ \{ H_{\nu } \} }
\def\Iseq{ \{ I_{\nu } \} }
\def\Jseq{ \{ J_{\nu } \} }
\def\Kseq{ \{ K_{\nu } \} }
\def\Lseq{ \{ L_{\nu } \} }
\def\Mseq{ \{ M_{\nu } \} }
\def\Nseq{ \{ N_{\nu } \} }
\def\Oseq{ \{ O_{\nu } \} }
\def\Pseq{ \{ P_{\nu } \} }
\def\Qseq{ \{ Q_{\nu } \} }
\def\Rseq{ \{ R_{\nu } \} }
\def\Sseq{ \{ S_{\nu } \} }
\def\Tseq{ \{ T_{\nu } \} }
\def\Useq{ \{ U_{\nu } \} }
\def\Vseq{ \{ V_{\nu } \} }
\def\Wseq{ \{ W_{\nu } \} }
\def\Xseq{ \{ X_{\nu } \} }
\def\Yseq{ \{ Y_{\nu } \} }
\def\Zseq{ \{ Z_{\nu } \} }

\def\alphaseq { \{ \alpha _{\nu } \} }
\def\betaseq { \{ \beta _{\nu } \} }
\def\gammaseq { \{ \gamma _{\nu } \} }
\def\deltaseq { \{ \delta _{\nu } \} }
\def\epsilonseq { \{ \epsilon _{\nu } \} }
\def\varepsilonseq { \{ \varepsilon _{\nu } \} }
\def\zetaseq { \{ \zeta _{\nu } \} }
\def\etaseq { \{ \eta _{\nu } \} }
\def\thetaseq { \{ \theta _{\nu } \} }
\def\varthetaseq { \{ \vartheta _{\nu } \} }
\def\iotaseq { \{ \iota _{\nu } \} }
\def\kappaseq { \{ \kappa _{\nu } \} }
\def\lambdaseq { \{ \lambda _{\nu } \} }
\def\museq { \{ \mu _{\nu } \} }
\def\nuseq { \{ \nu  _{\nu } \} }
\def\xiseq { \{ \xi _{\nu } \} }
\def\piseq { \{ \pi _{\nu } \} }
\def\varpiseq { \{ \varpi _{\nu } \} }
\def\rhoseq { \{ \rho _{\nu } \} }
\def\varrhoseq { \{ \varrho _{\nu } \} }
\def\sigmaseq { \{ \sigma _{\nu } \} }
\def\varsigmaseq { \{ \varsigma _{\nu } \} }
\def\tauseq { \{ \tau _{\nu } \} }
\def\upsilonseq { \{ \upsilon _{\nu } \} }
\def\phiseq { \{ \phi _{\nu } \} }
\def\varphiseq { \{ \varphi _{\nu } \} }
\def\chiseq { \{ \chi _{\nu } \} }
\def\psiseq { \{ \psi _{\nu } \} }
\def\omegaseq { \{ \omega _{\nu } \} }

\def\Gammaseq { \{ \Gamma _{\nu } \} }
\def\Deltaseq { \{ \Delta _{\nu } \} }
\def\Thetaseq { \{ \Theta _{\nu } \} }
\def\Lambdaseq { \{ \Lambda _{\nu } \} }
\def\Xiseq { \{ \Xi _{\nu } \} }
\def\Piseq { \{ \Pi _{\nu } \} }
\def\Sigmaseq { \{ \Sigma _{\nu } \} }
\def\Upsilonseq { \{ \Upsilon _{\nu } \} }
\def\Phiseq { \{ \Phi _{\nu } \} }
\def\Psiseq { \{ \Psi _{\nu } \} }
\def\Omegaseq { \{ \Omega _{\nu } \} }

\def\varGammaseq { \{ \varGamma _{\nu } \} }
\def\varDeltaseq { \{ \varDelta _{\nu } \} }
\def\varThetaseq { \{ \varTheta _{\nu } \} }
\def\varLambdaseq { \{ \varLambda _{\nu } \} }
\def\varXiseq { \{ \varXi _{\nu } \} }
\def\varPiseq { \{ \varPi _{\nu } \} }
\def\varSigmaseq { \{ \varSigma _{\nu } \} }
\def\varUpsilonseq { \{ \varUpsilon _{\nu } \} }
\def\varPhiseq { \{ \varPhi _{\nu } \} }
\def\varPsiseq { \{ \varPsi _{\nu } \} }
\def\varOmegaseq { \{ \varOmega _{\nu } \} }

\def\boldGammaseq { \{ \boldGamma _{\nu } \} }
\def\boldDeltaseq { \{ \boldDelta _{\nu } \} }
\def\boldThetaseq { \{ \boldTheta _{\nu } \} }
\def\boldLambdaseq { \{ \boldLambda _{\nu } \} }
\def\boldXiseq { \{ \boldXi _{\nu } \} }
\def\boldPiseq { \{ \boldPi _{\nu } \} }
\def\boldSigmaseq { \{ \boldSigma _{\nu } \} }
\def\boldUpsilonseq { \{ \boldUpsilon _{\nu } \} }
\def\boldPhiseq { \{ \boldPhi _{\nu } \} }
\def\boldPsiseq { \{ \boldPsi _{\nu } \} }
\def\boldOmegaseq { \{ \boldOmega _{\nu } \} }

\def\amu{   a_{\mu }  }
\def\bmu{   b_{\mu }  }
\def\cmu{   c_{\mu }  }
\def\dmu{   d_{\mu }  }
\def\emu{   e_{\mu }  }
\def\fmu{   f_{\mu }  }
\def\gmu{   g_{\mu }  }
\def\hmu{   h_{\mu }  }
\def\imu{   i_{\mu }  }
\def\jmu{   j_{\mu }  }
\def\kmu{   k_{\mu }  }
\def\lmu{   l_{\mu }  }
\def\mmu{   m_{\mu }  }
\def\nmu{   n_{\mu }  }
\def\omu{   o_{\mu }  }
\def\pmu{   p_{\mu }  }
\def\qmu{   q_{\mu }  }
\def\rmu{   r_{\mu }  }
\def\smu{   s_{\mu }  }
\def\tmu{   t_{\mu }  }
\def\umu{   u_{\mu }  }
\def\vmu{   v_{\mu }  }
\def\wmu{   w_{\mu }  }
\def\xmu{   x_{\mu }  }
\def\ymu{   y_{\mu }  }
\def\zmu{   z_{\mu }  }

\def\Amu{   A_{\mu }  }
\def\Bmu{   B_{\mu }  }
\def\Cmu{   C_{\mu }  }
\def\Dmu{   D_{\mu }  }
\def\Emu{   E_{\mu }  }
\def\Fmu{   F_{\mu }  }
\def\Gmu{   G_{\mu }  }
\def\Hmu{   H_{\mu }  }
\def\Imu{   I_{\mu }  }
\def\Jmu{   J_{\mu }  }
\def\Kmu{   K_{\mu }  }
\def\Lmu{   L_{\mu }  }
\def\Mmu{   M_{\mu }  }
\def\Nmu{   N_{\mu }  }
\def\Omu{   O_{\mu }  }
\def\Pmu{   P_{\mu }  }
\def\Qmu{   Q_{\mu }  }
\def\Rmu{   R_{\mu }  }
\def\Smu{   S_{\mu }  }
\def\Tmu{   T_{\mu }  }
\def\Umu{   U_{\mu }  }
\def\Vmu{   V_{\mu }  }
\def\Wmu{   W_{\mu }  }
\def\Xmu{   X_{\mu }  }
\def\Ymu{   Y_{\mu }  }
\def\Zmu{   Z_{\mu }  }

\def\alphamu{\alpha _{\mu }}
\def\betamu{\beta _{\mu }}
\def\gammamu{\gamma _{\mu }}
\def\deltamu{\delta _{\mu }}
\def\epsilonmu{\epsilon _{\mu }}
\def\varepsilonmu{\varepsilon _{\mu }}
\def\zetamu{\zeta _{\mu }}
\def\etamu{\eta _{\mu }}
\def\thetamu{\theta _{\mu }}
\def\varthetamu{\vartheta _{\mu }}
\def\iotamu{\iota _{\mu }}
\def\kappamu{\kappa _{\mu }}
\def\lambdamu{\lambda _{\mu }}
\def\mumu{\mu _{\mu }}
\def\numu{\nu _{\mu }}
\def\ximu{\xi _{\mu }}
\def\pimu{\pi _{\mu }}
\def\varpimu{\varpi _{\mu }}
\def\rhomu{\rho _{\mu }}
\def\varrhomu{\varrho _{\mu }}
\def\sigmamu{\sigma _{\mu }}
\def\varsigmamu{\varsigma _{\mu }}
\def\taumu{\tau _{\mu }}
\def\upsilonmu{\upsilon _{\mu }}
\def\phimu{\phi _{\mu }}
\def\varphimu{\varphi _{\mu }}
\def\chimu{\chi _{\mu }}
\def\psimu{\psi _{\mu }}
\def\omegamu{\omega _{\mu }}

\def\Gammamu{\Gamma _{\mu }}
\def\Deltamu{\Delta _{\mu }}
\def\Thetamu{\Theta _{\mu }}
\def\Lambdamu{\Lambda _{\mu }}
\def\Ximu{\Xi _{\mu }}
\def\Pimu{\Pi _{\mu }}
\def\Sigmamu{\Sigma _{\mu }}
\def\Upsilonmu{\Upsilon _{\mu }}
\def\Phimu{\Phi _{\mu }}
\def\Psimu{\Psi _{\mu }}
\def\Omegamu{\Omega _{\mu }}

\def\varGammamu{\varGamma _{\mu }}
\def\varDeltamu{\varDelta _{\mu }}
\def\varThetamu{\varTheta _{\mu }}
\def\varLambdamu{\varLambda _{\mu }}
\def\varXimu{\varXi _{\mu }}
\def\varPimu{\varPi _{\mu }}
\def\varSigmamu{\varSigma _{\mu }}
\def\varUpsilonmu{\varUpsilon _{\mu }}
\def\varPhimu{\varPhi _{\mu }}
\def\varPsimu{\varPsi _{\mu }}
\def\varOmegamu{\varOmega _{\mu }}

\def\boldGammamu{\boldGamma _{\mu }}
\def\boldDeltamu{\boldDelta _{\mu }}
\def\boldThetamu{\boldTheta _{\mu }}
\def\boldLambdamu{\boldLambda _{\mu }}
\def\boldXimu{\boldXi _{\mu }}
\def\boldPimu{\boldPi _{\mu }}
\def\boldSigmamu{\boldSigma _{\mu }}
\def\boldUpsilonmu{\boldUpsilon _{\mu }}
\def\boldPhimu{\boldPhi _{\mu }}
\def\boldPsimu{\boldPsi _{\mu }}
\def\boldOmegamu{\boldOmega _{\mu }}


\def\asmu{   a^{(\mu )}  }
\def\bsmu{   b^{(\mu )}  }
\def\csmu{   c^{(\mu )}  }
\def\dsmu{   d^{(\mu )}  }
\def\esmu{   e^{(\mu )}  }
\def\fsmu{   f^{(\mu )}  }
\def\gsmu{   g^{(\mu )}  }
\def\hsmu{   h^{(\mu )}  }
\def\ismu{   i^{(\mu )}  }
\def\jsmu{   j^{(\mu )}  }
\def\ksmu{   k^{(\mu )}  }
\def\lsmu{   l^{(\mu )}  }
\def\msmu{   m^{(\mu )}  }
\def\nsmu{   n^{(\mu )}  }
\def\osmu{   o^{(\mu )}  }
\def\psmu{   p^{(\mu )}  }
\def\qsmu{   q^{(\mu )}  }
\def\rsmu{   r^{(\mu )}  }
\def\ssmu{   s^{(\mu )}  }
\def\tsmu{   t^{(\mu )}  }
\def\usmu{   u^{(\mu )}  }
\def\vsmu{   v^{(\mu )}  }
\def\wsmu{   w^{(\mu )}  }
\def\xsmu{   x^{(\mu )}  }
\def\ysmu{   y^{(\mu )}  }
\def\zsmu{   z^{(\mu )}  }

\def\Asmu{   A^{(\mu )}  }
\def\Bsmu{   B^{(\mu )}  }
\def\Csmu{   C^{(\mu )}  }
\def\Dsmu{   D^{(\mu )}  }
\def\Esmu{   E^{(\mu )}  }
\def\Fsmu{   F^{(\mu )}  }
\def\Gsmu{   G^{(\mu )}  }
\def\Hsmu{   H^{(\mu )}  }
\def\Ismu{   I^{(\mu )}  }
\def\Jsmu{   J^{(\mu )}  }
\def\Ksmu{   K^{(\mu )}  }
\def\Lsmu{   L^{(\mu )}  }
\def\Msmu{   M^{(\mu )}  }
\def\Nsmu{   N^{(\mu )}  }
\def\Osmu{   O^{(\mu )}  }
\def\Psmu{   P^{(\mu )}  }
\def\Qsmu{   Q^{(\mu )}  }
\def\Rsmu{   R^{(\mu )}  }
\def\Ssmu{   S^{(\mu )}  }
\def\Tsmu{   T^{(\mu )}  }
\def\Usmu{   U^{(\mu )}  }
\def\Vsmu{   V^{(\mu )}  }
\def\Wsmu{   W^{(\mu )}  }
\def\Xsmu{   X^{(\mu )}  }
\def\Ysmu{   Y^{(\mu )}  }
\def\Zsmu{   Z^{(\mu )}  }

\def\alphasmu{\alpha ^{(\mu )}}
\def\betasmu{\beta ^{(\mu )}}
\def\gammasmu{\gamma ^{(\mu )}}
\def\deltasmu{\delta ^{(\mu )}}
\def\epsilonsmu{\epsilon ^{(\mu )}}
\def\varepsilonsmu{\varepsilon ^{(\mu )}}
\def\zetasmu{\zeta ^{(\mu )}}
\def\etasmu{\eta ^{(\mu )}}
\def\thetasmu{\theta ^{(\mu )}}
\def\varthetasmu{\vartheta ^{(\mu )}}
\def\iotasmu{\iota ^{(\mu )}}
\def\kappasmu{\kappa ^{(\mu )}}
\def\lambdasmu{\lambda ^{(\mu )}}
\def\musmu{\mu ^{(\mu )}}
\def\nusmu{\nu ^{(\mu )}}
\def\xismu{\xi ^{(\mu )}}
\def\pismu{\pi ^{(\mu )}}
\def\varpismu{\varpi ^{(\mu )}}
\def\rhosmu{\rho ^{(\mu )}}
\def\varrhosmu{\varrho ^{(\mu )}}
\def\sigmasmu{\sigma ^{(\mu )}}
\def\varsigmasmu{\varsigma ^{(\mu )}}
\def\tausmu{\tau ^{(\mu )}}
\def\upsilonsmu{\upsilon ^{(\mu )}}
\def\phismu{\phi ^{(\mu )}}
\def\varphismu{\varphi ^{(\mu )}}
\def\chismu{\chi ^{(\mu )}}
\def\psismu{\psi ^{(\mu )}}
\def\omegasmu{\omega ^{(\mu )}}

\def\Gammasmu{\Gamma ^{(\mu )}}
\def\Deltasmu{\Delta ^{(\mu )}}
\def\Thetasmu{\Theta ^{(\mu )}}
\def\Lambdasmu{\Lambda ^{(\mu )}}
\def\Xismu{\Xi ^{(\mu )}}
\def\Pismu{\Pi ^{(\mu )}}
\def\Sigmasmu{\Sigma ^{(\mu )}}
\def\Upsilonsmu{\Upsilon ^{(\mu )}}
\def\Phismu{\Phi ^{(\mu )}}
\def\Psismu{\Psi ^{(\mu )}}
\def\Omegasmu{\Omega ^{(\mu )}}

\def\varGammasmu{\varGamma ^{(\mu )}}
\def\varDeltasmu{\varDelta ^{(\mu )}}
\def\varThetasmu{\varTheta ^{(\mu )}}
\def\varLambdasmu{\varLambda ^{(\mu )}}
\def\varXismu{\varXi ^{(\mu )}}
\def\varPismu{\varPi ^{(\mu )}}
\def\varSigmasmu{\varSigma ^{(\mu )}}
\def\varUpsilonsmu{\varUpsilon ^{(\mu )}}
\def\varPhismu{\varPhi ^{(\mu )}}
\def\varPsismu{\varPsi ^{(\mu )}}
\def\varOmegasmu{\varOmega ^{(\mu )}}

\def\boldGammasmu{\boldGamma ^{(\mu )}}
\def\boldDeltasmu{\boldDelta ^{(\mu )}}
\def\boldThetasmu{\boldTheta ^{(\mu )}}
\def\boldLambdasmu{\boldLambda ^{(\mu )}}
\def\boldXismu{\boldXi ^{(\mu )}}
\def\boldPismu{\boldPi ^{(\mu )}}
\def\boldSigmasmu{\boldSigma ^{(\mu )}}
\def\boldUpsilonsmu{\boldUpsilon ^{(\mu )}}
\def\boldPhismu{\boldPhi ^{(\mu )}}
\def\boldPsismu{\boldPsi ^{(\mu )}}
\def\boldOmegasmu{\boldOmega ^{(\mu )}}

\def\anu{   a_{\nu }  }
\def\bnu{   b_{\nu }  }
\def\cnu{   c_{\nu }  }
\def\dnu{   d_{\nu }  }
\def\enu{   e_{\nu }  }
\def\fnu{   f_{\nu }  }
\def\gnu{   g_{\nu }  }
\def\hnu{   h_{\nu }  }
\def\inu{   i_{\nu }  }
\def\jnu{   j_{\nu }  }
\def\knu{   k_{\nu }  }
\def\lnu{   l_{\nu }  }
\def\mnu{   m_{\nu }  }
\def\nnu{   n_{\nu }  }
\def\onu{   o_{\nu }  }
\def\pnu{   p_{\nu }  }
\def\qnu{   q_{\nu }  }
\def\rnu{   r_{\nu }  }
\def\snu{   s_{\nu }  }
\def\tnu{   t_{\nu }  }
\def\unu{   u_{\nu }  }
\def\vnu{   v_{\nu }  }
\def\wnu{   w_{\nu }  }
\def\xnu{   x_{\nu }  }
\def\ynu{   y_{\nu }  }
\def\znu{   z_{\nu }  }

\def\Anu{   A_{\nu }  }
\def\Bnu{   B_{\nu }  }
\def\Cnu{   C_{\nu }  }
\def\Dnu{   D_{\nu }  }
\def\Enu{   E_{\nu }  }
\def\Fnu{   F_{\nu }  }
\def\Gnu{   G_{\nu }  }
\def\Hnu{   H_{\nu }  }
\def\Inu{   I_{\nu }  }
\def\Jnu{   J_{\nu }  }
\def\Knu{   K_{\nu }  }
\def\Lnu{   L_{\nu }  }
\def\Mnu{   M_{\nu }  }
\def\Nnu{   N_{\nu }  }
\def\Onu{   O_{\nu }  }
\def\Pnu{   P_{\nu }  }
\def\Qnu{   Q_{\nu }  }
\def\Rnu{   R_{\nu }  }
\def\Snu{   S_{\nu }  }
\def\Tnu{   T_{\nu }  }
\def\Unu{   U_{\nu }  }
\def\Vnu{   V_{\nu }  }
\def\Wnu{   W_{\nu }  }
\def\Xnu{   X_{\nu }  }
\def\Ynu{   Y_{\nu }  }
\def\Znu{   Z_{\nu }  }

\def\alphanu{\alpha _{\nu }}
\def\betanu{\beta _{\nu }}
\def\gammanu{\gamma _{\nu }}
\def\deltanu{\delta _{\nu }}
\def\epsilonnu{\epsilon _{\nu }}
\def\varepsilonnu{\varepsilon _{\nu }}
\def\zetanu{\zeta _{\nu }}
\def\etanu{\eta _{\nu }}
\def\thetanu{\theta _{\nu }}
\def\varthetanu{\vartheta _{\nu }}
\def\iotanu{\iota _{\nu }}
\def\kappanu{\kappa _{\nu }}
\def\lambdanu{\lambda _{\nu }}
\def\munu{\mu _{\nu }}
\def\nunu{\nu _{\nu }}
\def\xinu{\xi _{\nu }}
\def\pinu{\pi _{\nu }}
\def\varpinu{\varpi _{\nu }}
\def\rhonu{\rho _{\nu }}
\def\varrhonu{\varrho _{\nu }}
\def\sigmanu{\sigma _{\nu }}
\def\varsigmanu{\varsigma _{\nu }}
\def\taunu{\tau _{\nu }}
\def\upsilonnu{\upsilon _{\nu }}
\def\phinu{\phi _{\nu }}
\def\varphinu{\varphi _{\nu }}
\def\chinu{\chi _{\nu }}
\def\psinu{\psi _{\nu }}
\def\omeganu{\omega _{\nu }}

\def\Gammanu{\Gamma _{\nu }}
\def\Deltanu{\Delta _{\nu }}
\def\Thetanu{\Theta _{\nu }}
\def\Lambdanu{\Lambda _{\nu }}
\def\Xinu{\Xi _{\nu }}
\def\Pinu{\Pi _{\nu }}
\def\Sigmanu{\Sigma _{\nu }}
\def\Upsilonnu{\Upsilon _{\nu }}
\def\Phinu{\Phi _{\nu }}
\def\Psinu{\Psi _{\nu }}
\def\Omeganu{\Omega _{\nu }}

\def\varGammanu{\varGamma _{\nu }}
\def\varDeltanu{\varDelta _{\nu }}
\def\varThetanu{\varTheta _{\nu }}
\def\varLambdanu{\varLambda _{\nu }}
\def\varXinu{\varXi _{\nu }}
\def\varPinu{\varPi _{\nu }}
\def\varSigmanu{\varSigma _{\nu }}
\def\varUpsilonnu{\varUpsilon _{\nu }}
\def\varPhinu{\varPhi _{\nu }}
\def\varPsinu{\varPsi _{\nu }}
\def\varOmeganu{\varOmega _{\nu }}

\def\boldGammanu{\boldGamma _{\nu }}
\def\boldDeltanu{\boldDelta _{\nu }}
\def\boldThetanu{\boldTheta _{\nu }}
\def\boldLambdanu{\boldLambda _{\nu }}
\def\boldXinu{\boldXi _{\nu }}
\def\boldPinu{\boldPi _{\nu }}
\def\boldSigmanu{\boldSigma _{\nu }}
\def\boldUpsilonnu{\boldUpsilon _{\nu }}
\def\boldPhinu{\boldPhi _{\nu }}
\def\boldPsinu{\boldPsi _{\nu }}
\def\boldOmeganu{\boldOmega _{\nu }}


\def\asnu{   a^{(\nu )}  }
\def\bsnu{   b^{(\nu )}  }
\def\csnu{   c^{(\nu )}  }
\def\dsnu{   d^{(\nu )}  }
\def\esnu{   e^{(\nu )}  }
\def\fsnu{   f^{(\nu )}  }
\def\gsnu{   g^{(\nu )}  }
\def\hsnu{   h^{(\nu )}  }
\def\isnu{   i^{(\nu )}  }
\def\jsnu{   j^{(\nu )}  }
\def\ksnu{   k^{(\nu )}  }
\def\lsnu{   l^{(\nu )}  }
\def\msnu{   m^{(\nu )}  }
\def\nsnu{   n^{(\nu )}  }
\def\osnu{   o^{(\nu )}  }
\def\psnu{   p^{(\nu )}  }
\def\qsnu{   q^{(\nu )}  }
\def\rsnu{   r^{(\nu )}  }
\def\ssnu{   s^{(\nu )}  }
\def\tsnu{   t^{(\nu )}  }
\def\usnu{   u^{(\nu )}  }
\def\vsnu{   v^{(\nu )}  }
\def\wsnu{   w^{(\nu )}  }
\def\xsnu{   x^{(\nu )}  }
\def\ysnu{   y^{(\nu )}  }
\def\zsnu{   z^{(\nu )}  }

\def\Asnu{   A^{(\nu )}  }
\def\Bsnu{   B^{(\nu )}  }
\def\Csnu{   C^{(\nu )}  }
\def\Dsnu{   D^{(\nu )}  }
\def\Esnu{   E^{(\nu )}  }
\def\Fsnu{   F^{(\nu )}  }
\def\Gsnu{   G^{(\nu )}  }
\def\Hsnu{   H^{(\nu )}  }
\def\Isnu{   I^{(\nu )}  }
\def\Jsnu{   J^{(\nu )}  }
\def\Ksnu{   K^{(\nu )}  }
\def\Lsnu{   L^{(\nu )}  }
\def\Msnu{   M^{(\nu )}  }
\def\Nsnu{   N^{(\nu )}  }
\def\Osnu{   O^{(\nu )}  }
\def\Psnu{   P^{(\nu )}  }
\def\Qsnu{   Q^{(\nu )}  }
\def\Rsnu{   R^{(\nu )}  }
\def\Ssnu{   S^{(\nu )}  }
\def\Tsnu{   T^{(\nu )}  }
\def\Usnu{   U^{(\nu )}  }
\def\Vsnu{   V^{(\nu )}  }
\def\Wsnu{   W^{(\nu )}  }
\def\Xsnu{   X^{(\nu )}  }
\def\Ysnu{   Y^{(\nu )}  }
\def\Zsnu{   Z^{(\nu )}  }

\def\alphasnu{\alpha ^{(\nu )}}
\def\betasnu{\beta ^{(\nu )}}
\def\gammasnu{\gamma ^{(\nu )}}
\def\deltasnu{\delta ^{(\nu )}}
\def\epsilonsnu{\epsilon ^{(\nu )}}
\def\varepsilonsnu{\varepsilon ^{(\nu )}}
\def\zetasnu{\zeta ^{(\nu )}}
\def\etasnu{\eta ^{(\nu )}}
\def\thetasnu{\theta ^{(\nu )}}
\def\varthetasnu{\vartheta ^{(\nu )}}
\def\iotasnu{\iota ^{(\nu )}}
\def\kappasnu{\kappa ^{(\nu )}}
\def\lambdasnu{\lambda ^{(\nu )}}
\def\musnu{\mu ^{(\nu )}}
\def\nusnu{\nu ^{(\nu )}}
\def\xisnu{\xi ^{(\nu )}}
\def\pisnu{\pi ^{(\nu )}}
\def\varpisnu{\varpi ^{(\nu )}}
\def\rhosnu{\rho ^{(\nu )}}
\def\varrhosnu{\varrho ^{(\nu )}}
\def\sigmasnu{\sigma ^{(\nu )}}
\def\varsigmasnu{\varsigma ^{(\nu )}}
\def\tausnu{\tau ^{(\nu )}}
\def\upsilonsnu{\upsilon ^{(\nu )}}
\def\phisnu{\phi ^{(\nu )}}
\def\varphisnu{\varphi ^{(\nu )}}
\def\chisnu{\chi ^{(\nu )}}
\def\psisnu{\psi ^{(\nu )}}
\def\omegasnu{\omega ^{(\nu )}}

\def\Gammasnu{\Gamma ^{(\nu )}}
\def\Deltasnu{\Delta ^{(\nu )}}
\def\Thetasnu{\Theta ^{(\nu )}}
\def\Lambdasnu{\Lambda ^{(\nu )}}
\def\Xisnu{\Xi ^{(\nu )}}
\def\Pisnu{\Pi ^{(\nu )}}
\def\Sigmasnu{\Sigma ^{(\nu )}}
\def\Upsilonsnu{\Upsilon ^{(\nu )}}
\def\Phisnu{\Phi ^{(\nu )}}
\def\Psisnu{\Psi ^{(\nu )}}
\def\Omegasnu{\Omega ^{(\nu )}}

\def\varGammasnu{\varGamma ^{(\nu )}}
\def\varDeltasnu{\varDelta ^{(\nu )}}
\def\varThetasnu{\varTheta ^{(\nu )}}
\def\varLambdasnu{\varLambda ^{(\nu )}}
\def\varXisnu{\varXi ^{(\nu )}}
\def\varPisnu{\varPi ^{(\nu )}}
\def\varSigmasnu{\varSigma ^{(\nu )}}
\def\varUpsilonsnu{\varUpsilon ^{(\nu )}}
\def\varPhisnu{\varPhi ^{(\nu )}}
\def\varPsisnu{\varPsi ^{(\nu )}}
\def\varOmegasnu{\varOmega ^{(\nu )}}

\def\boldGammasnu{\boldGamma ^{(\nu )}}
\def\boldDeltasnu{\boldDelta ^{(\nu )}}
\def\boldThetasnu{\boldTheta ^{(\nu )}}
\def\boldLambdasnu{\boldLambda ^{(\nu )}}
\def\boldXisnu{\boldXi ^{(\nu )}}
\def\boldPisnu{\boldPi ^{(\nu )}}
\def\boldSigmasnu{\boldSigma ^{(\nu )}}
\def\boldUpsilonsnu{\boldUpsilon ^{(\nu )}}
\def\boldPhisnu{\boldPhi ^{(\nu )}}
\def\boldPsisnu{\boldPsi ^{(\nu )}}
\def\boldOmegasnu{\boldOmega ^{(\nu )}}


\def\vphi {\varphi }


\def\inv{   ^{-1}  }

\def\ainv{   a^{-1}  }
\def\binv{   b^{-1}  }
\def\cinv{   c^{-1}  }
\def\dinv{   d^{-1}  }
\def\einv{   e^{-1}  }
\def\finv{   f^{-1}  }
\def\ginv{   g^{-1}  }
\def\hinv{   h^{-1}  }
\def\iinv{   i^{-1}  }
\def\jinv{   j^{-1}  }
\def\kinv{   k^{-1}  }
\def\linv{   l^{-1}  }
\def\minv{   m^{-1}  }
\def\ninv{   n^{-1}  }
\def\oinv{   o^{-1}  }
\def\pinv{   p^{-1}  }
\def\qinv{   q^{-1}  }
\def\rinv{   r^{-1}  }
\def\sinv{   s^{-1}  }
\def\tinv{   t^{-1}  }
\def\uinv{   u^{-1}  }
\def\vinv{   v^{-1}  }
\def\winv{   w^{-1}  }
\def\xinv{   x^{-1}  }
\def\yinv{   y^{-1}  }
\def\zinv{   z^{-1}  }

\def\Ainv{   A^{-1}  }
\def\Binv{   B^{-1}  }
\def\Cinv{   C^{-1}  }
\def\Dinv{   D^{-1}  }
\def\Einv{   E^{-1}  }


\def\Ginv{   G^{-1}  }
\def\Hinv{   H^{-1}  }
\def\Iinv{   I^{-1}  }
\def\Jinv{   J^{-1}  }
\def\Kinv{   K^{-1}  }
\def\Linv{   L^{-1}  }
\def\Minv{   M^{-1}  }
\def\Ninv{   N^{-1}  }
\def\Oinv{   O^{-1}  }
\def\Pinv{   P^{-1}  }
\def\Qinv{   Q^{-1}  }
\def\Rinv{   R^{-1}  }
\def\Sinv{   S^{-1}  }
\def\Tinv{   T^{-1}  }
\def\Uinv{   U^{-1}  }
\def\Vinv{   V^{-1}  }
\def\Winv{   W^{-1}  }
\def\Xinv{   X^{-1}  }
\def\Yinv{   Y^{-1}  }
\def\Zinv{   Z^{-1}  }

\def\alphainv{\alpha ^{-1}}
\def\betainv{\beta ^{-1}}
\def\gammainv{\gamma ^{-1}}
\def\deltainv{\delta ^{-1}}
\def\epsiloninv{\epsilon ^{-1}}
\def\varepsiloninv{\varepsilon ^{-1}}
\def\zetainv{\zeta ^{-1}}
\def\etainv{\eta ^{-1}}
\def\thetainv{\theta ^{-1}}
\def\varthetainv{\vartheta ^{-1}}
\def\iotainv{\iota ^{-1}}
\def\kappainv{\kappa ^{-1}}
\def\lambdainv{\lambda ^{-1}}
\def\muinv{\mu ^{-1}}
\def\nuinv{\nu ^{-1}}
\def\xiinv{\xi ^{-1}}
\def\piinv{\pi ^{-1}}
\def\varpiinv{\varpi ^{-1}}
\def\rhoinv{\rho ^{-1}}
\def\varrhoinv{\varrho ^{-1}}
\def\sigmainv{\sigma ^{-1}}
\def\varsigmainv{\varsigma ^{-1}}
\def\tauinv{\tau ^{-1}}
\def\upsiloninv{\upsilon ^{-1}}
\def\phiinv{\phi ^{-1}}
\def\varphiinv{\varphi ^{-1}}
\def\vphiinv{\varphi ^{-1}}
\def\chiinv{\chi ^{-1}}
\def\psiinv{\psi ^{-1}}
\def\omegainv{\omega ^{-1}}

\def\Gammainv{\Gamma ^{-1}}
\def\Deltainv{\Delta ^{-1}}
\def\Thetainv{\Theta ^{-1}}
\def\Lambdainv{\Lambda ^{-1}}
\def\Xiinv{\Xi ^{-1}}
\def\Piinv{\Pi ^{-1}}
\def\Sigmainv{\Sigma ^{-1}}
\def\Upsiloninv{\Upsilon ^{-1}}
\def\Phiinv{\Phi ^{-1}}
\def\Psiinv{\Psi ^{-1}}
\def\Omegainv{\Omega ^{-1}}

\def\varGammainv{\varGamma ^{-1}}
\def\varDeltainv{\varDelta ^{-1}}
\def\varThetainv{\varTheta ^{-1}}
\def\varLambdainv{\varLambda ^{-1}}
\def\varXiinv{\varXi ^{-1}}
\def\varPiinv{\varPi ^{-1}}
\def\varSigmainv{\varSigma ^{-1}}
\def\varUpsiloninv{\varUpsilon ^{-1}}
\def\varPhiinv{\varPhi ^{-1}}
\def\varPsiinv{\varPsi ^{-1}}
\def\varOmegainv{\varOmega ^{-1}}

\def\boldGammainv{\boldGamma ^{-1}}
\def\boldDeltainv{\boldDelta ^{-1}}
\def\boldThetainv{\boldTheta ^{-1}}
\def\boldLambdainv{\boldLambda ^{-1}}
\def\boldXiinv{\boldXi ^{-1}}
\def\boldPiinv{\boldPi ^{-1}}
\def\boldSigmainv{\boldSigma ^{-1}}
\def\boldUpsiloninv{\boldUpsilon ^{-1}}
\def\boldPhiinv{\boldPhi ^{-1}}
\def\boldPsiinv{\boldPsi ^{-1}}
\def\boldOmegainv{\boldOmega ^{-1}}
\title[Thompson's group $F$ is not K\"ahler]
{Thompson's group $F$ is not K\"ahler}
\author[T.~Napier]{Terrence Napier$^{*}$}
\address{Department of Mathematics\\Lehigh University\\Bethlehem, PA 18015}
\email{tjn2@lehigh.edu}
\thanks{$^{*}$Research partially
supported by NSF grant DMS0306441}
\author[M.~Ramachandran]{Mohan Ramachandran}
\address{Department of Mathematics\\SUNY at Buffalo\\Buffalo, NY 14260}
\email{ramac-m@newton.math.buffalo.edu}



\dedicatory{Dedicated to Ross Geoghegan in honor of his
60$^{\text{th}}$ birthday.}

\date{May 25, 2005}

\begin{abstract}
The purpose of this note is to prove that Richard Thompson's group
$F$ and variants of it studied by Ken Brown are not K\"ahler
groups.
\end{abstract}

\maketitle

\section*{Introduction} \label{introduction}

The purpose of this note is to prove the following:
\begin{thm} Thompson's group $F$ and the generalizations
$F_{n,\infty }$ and $F_n$ for $n=2,3,4, \dots $ are not K\"ahler.
\end{thm}
Theorem~0.1 answers a question of Ross Geoghegan (see Ken Brown's
paper in these proceedings \cite{Br2}).

Thompson discovered the group $F$ in 1965 in the context of his
work in mathematical logic. Brown and Geoghegan \cite{BrGe}
determined that $F$ is of type $FP_\infty$, thus making $F$ the
first known example of a torsion free group of type $FP_\infty$
which is not of type $FP$. The group $F$ has also appeared in
homotopy theory \cite{FH}. For more details on the properties and
history of Thompson's groups $F$, $F_{n,\infty }$, and $F_n$, the
reader may refer to \cite{CFP}, \cite{Br1}, \cite{BG}, and
\cite{BrGe}.

A finitely presented group is called a {\it K\"ahler group} if it
is the fundamental group of a \cpt K\"ahler manifold.  A central
problem in the study of the topology of \cpt K\"ahler manifolds
(for example, smooth projective varieties) is that of determining
which groups are K\"ahler groups.  For example, according to Hodge
theory, the Abelianization of a K\"ahler group must be of even
rank (see, for example, \cite{W}). By \cite{ArBR}, a K\"ahler
group has at most one end. If $M$ is any Hopf surface, then $M$ is
a non-K\"ahler \cpt \cpx manifold  which has fundamental
group~$\Z$ and $2$-ended universal covering $\C^2\setminus\{ 0\}$.
The Heisenberg group $H$ of $3\times 3$ upper triangular integer
matrices with diagonal entries~$1$ has exactly one end and its
Abelianization is of rank~$2$, but $H$ is not a K\"ahler group
(see, for example, \cite{JR} or \cite{Ar}). On the positive side,
any finite group is the fundamental group of a smooth projective
variety. For any positive integer~$g$, the group $\langle
a_1,a_2,\dots , a_{2g}\mid [a_1,a_{g+1}]\cdot\cdots\cdot
[a_g,a_{2g}]=1\rangle $ is the fundamental group of a curve of
genus~$g$. In fact, much more subtle examples are now known to
exist; for example, Toledo's examples of K\"ahler groups which are
not residually finite~\cite{T}. The problem is, of course, related
to that of determining which groups are fundamental groups of
smooth projective varieties.  In fact, there are no known examples
of K\"ahler groups which are not also fundamental groups of smooth
projective varieties. For more details on the study of K\"ahler
groups, the reader may refer to the surveys \cite{AmBCKT} and
\cite{Ar}.

We give two proofs that Thompson's group $F$ is not K\"ahler. The
first proof (Section~1) relies on the particular properties
of~$F$. The second proof (Section~2) gives the following more
general fact which may eventually yield a proof that other groups
of interest are not K\"ahler:
\begin{thm}
Let $G$ be a group which satisfies the following:
\begin{enumerate}
\item [(i)] Any nontrivial normal subgroup contains the commutator
subgroup $G_{\text{{\rm com}}}$ (i.e. every proper quotient of $G$
is Abelian), and \item [(ii)] $G$ is a properly ascending HNN
extension.
\end{enumerate}
Then $G$ and any group containing $G$ as a subgroup of finite
index are not K\"ahler.
\end{thm}

Since $F_{n,\infty }$ ($F=F_{2,\infty }=F_2$) satisfies the above
conditions (i)~and~(ii) and $F_{n,\infty }$ is a (normal) subgroup
of finite index in $F_n$ for every $n\geq 2$ (see, for example,
\cite{CFP} and \cite{BG}), Theorem~0.1 is a consequence of
Theorem~0.2.

\begin{ack} We would like to thank Ken Brown and Ross Geoghegan
for answering our questions about Thompson's group $F$ and Ross
Geoghegan for asking the question that this paper answers. We
would also like to thank John Meier for trying to teach us right
from wrong in geometric group theory. Finally, we would like to
thank the referee for helpful comments.
\end{ack}

\section{First proof}

For the purposes of this note, $X$ will denote a \con \cpt
K\"ahler manifold and $C$ will denote a \con \cpt curve (i.e.~a
\cpt $1$-dimensional \cpx manifold). A remark we will use without
comment is a nonconstant \holo map from $X$ to $C$ is surjective
and open.

Let $f:X\to C$ be a surjective \holo map with \con fibers. Let $\{
\, p_1,\dots ,p_r\} $ be the set of critical values of $f$ and let
$m_i$ be the greatest common divisor of the multiplicities of the
\comps of the divisor $f\inv (p_i)$ for each $i=1,\dots , r$.  Let
$C^{\text{orb}}$ be the $2$-orbifold with underlying topological
space $C$ and singular points $p_1,\dots , p_r$ of order
$m_1,\dots ,m_r$, respectively. Note that $C^{\text{orb}}$ has a
unique structure of a \cpx $1$-orbifold so that the map $f$ is
\holons . We recall that, in this situation (see \cite{Ct},
\cite{Si}), the {\it orbifold fundamental group} $\pi
_1^{\text{orb}}(C)$ is the quotient of $\pi _1(C- \{ \, p_1,\dots
,p_r\, \} )$ by the normal closure of $\setof{ \gamma
_i^{m_i}}{i=1,\dots ,r}$, where $\gamma _i$ is a simple loop
around the point $p_i$ for each $i$.

Let $G=\pi_1(X)$. The kernel $K$ of the surjective homomorphism
$\pi_1(X)\to \pi _1^{\text{orb}}(C)$ is the image in $G$ of the
fundamental group of a general fiber of $f$. In particular, $K$ is
a finitely generated normal subgroup of $G$. If the underlying
topological surface of $C$ has positive genus, then $\pi
_1^{\text{orb}}(C)$ is $\Z^2$ or a co\cpt Fuchsian group. Summing
up, we have the following:
\begin{lem}
For $G=\pi _1(X)$, we have an exact sequence
$$
1\to K\to G\to \pi _1^{\text{orb}}(C)\to 1,
$$
where $K$ is finitely generated. If $C$ has positive genus, then
$\pi _1^{\text{orb}}(C)$ is either $\Z^2$ or a co\cpt Fuchsian
group.
\end{lem}

In both proofs of Theorem~0.1, the main point leading to a
contradiction is that, if $F=\pi _1(X)$ for $X$ K\"ahler, then $X$
admits such a mapping $f:X\to C$. In the first proof, the main
features of the group $F$ which we will use are (see [CFP]):
\begin{enumerate} \item [1.] $F$ is torsion free;
\item [2.] The commutator subgroup $F_{\text{com}}$ is {\it not}
finitely generated; \item [3.] Any nontrivial normal subgroup
contains $F_{\text{com}}$; and \item [4.] The Abelianization
$F/F_{\text{com}}$ is $\Z^2$.
\end{enumerate}

\begin{pf*}{First proof of Theorem~0.1}
If $F=\pi_1(X)$ for some \con \cpt K\"ahler manifold $X$, then,
since the Abelianization $F/F_{\text{com}}$ of $F$ is $\Z^2$
(property~4), the Albanese variety is an elliptic curve $E$. Stein
factoring the Albanese map $h:X\to E$, we get a commutative
diagram of surjective \holo maps
\begin{center}\begin{picture}(230,52)
 \put(128,0){$C$}
\put(80,39){$X$} \put(105,45){$h$} \put(88,43){\vector(1,0){37}}
\put(128,39){$E$} \put(136,21){$g$} \put(131,10){\vector(0,1){26}}
\put(99,14){$f$} \put(90,40){\vector(1,-1){32}}
\end{picture}
\end{center}
where $C$ is a \cpt curve, $f$ has \con fibers, and $g$ has finite
fibers. In particular, $C$ is of positive genus. Applying
Lemma~1.1, we get the exact sequence
$$
1\to K\to F=\pi_1(X)\to \pi _1^{\text{orb}}(C)\to 1,
$$
where $K$ is finitely generated. Thus $F/K\cong \pi
_1^{\text{orb}}(C)$.  Now $F$ is not co\cpt Fuchsian. For, if $F$
were co\cpt Fuchsian, then, since $F$ is torsion-free
(property~1), $F$ would be the fundamental group of a curve of
genus at least~$2$ and hence $F$ would have Abelianization of rank
at least~$4$ (contradicting property~4). $F$ also cannot be
isomorphic to $\Z^2$ because the subgroup $F_{\text{com}}$ is not
finitely generated (property~2). Thus $K$ is nontrivial and,
therefore, $K$ contains the commutator subgroup~$F_{\text{com}}$
(property~3). Therefore $\pi _1^{\text{orb}}(C)$ is a quotient of
the Abelianization $F/F_{\text{com}}\cong\Z^2$. In particular,
since co\cpt Fuchsian groups are non-Abelian, $\pi
_1^{\text{orb}}(C)$ must be isomorphic to $\Z^2$ and hence cannot
be isomorphic to a {\it proper} quotient of $\Z^2$.
%
%
Thus we must have $\pi _1^{\text{orb}}(C)=F/F_{\text{com}}$; that
is, $K=F_{\text{com}}$. But this is impossible because
$F_{\text{com}}$ is not finitely generated (property~2).
\end{pf*}

\section{Second proof}

We first consider the following:
\begin{thm}
Given a \con \cpt K\"ahler manifold~$X$ and an exact sequence
$$
(*)\qquad 1\to N\to \pi_1(X)\overset{\rho}{\to }\Z\to 1
$$
with $N$ not finitely generated, we get an exact sequence
$$
(**)\qquad 1\to K\to \pi_1(X)\to \Gamma\to 1,
$$
where $\Gamma $ is a co\cpt Fuchsian group of a curve of positive
genus and $K$ is finitely generated.
\end{thm}
\begin{pf}
By Theorem~4.3 of~\cite{NR}, there is a surjective \holo map with
\con fibers to a curve of positive genus $f:X\to C$ and a
factorization
\begin{center}\begin{picture}(230,64)
 \put(137,0){$\Z$}
\put(60,46){$\pi_1(X)$} \put(110,53){$f_*$}
\put(92,50){\vector(1,0){41}} \put(137,46){$\pi
_1^{\text{orb}}(C)$} \put(146,27){$\rho '$}
\put(140,42){\vector(0,-1){29}} \put(103,17){$\rho$}
\put(92,45){\vector(1,-1){40}}
\end{picture}
\end{center}
By Lemma~1.1, we have an exact sequence $1\to K\to \pi_1(X)\to \pi
_1^{\text{orb}}(C)\to 1$, where $K$ is finitely generated. Since
$K$ is finitely generated but $N$ is not, $\ker (\rho ')\cong N/K$
is not finitely generated. It follows that $\Gamma =\pi
_1^{\text{orb}}(C)$ is a co\cpt Fuchsian group, since every
subgroup of $\Z^2$ is finitely generated.
\end{pf}

\begin{pf*}{Proof of Theorem~0.2}
It suffices to prove that $G$ is not K\"ahler, because any finite
covering of a \cpt K\"ahler manifold is \cpt K\"ahler. If there
exists a \con \cpt K\"ahler manifold $X$ with $\pi_1(X)=G$, then,
since $G$ is a properly ascending HNN extension (property~(ii)),
we get an exact sequence of the form $(*)$ as in Theorem~2.1 and,
therefore, we get an exact sequence of the form $(**)$. If $K$ is
nontrivial, then $K$ contains the commutator subgroup
(property~(i)). But a co\cpt Fuchsian group cannot be Abelian, so
we arrive at a contradiction. If $K$ is trivial, then $G$ is
co\cpt Fuchsian and, therefore, by Malcev's theorem, residually
finite; i.e.~for any element $g\neq 1$, there is a finite
homomorphic image of $G$ in which the image of $g$ is not the
identity (for an elementary proof, see \cite{Al}). Taking $g\in
G_{\text{com}}\setminus\{ 1\}$, we get a non-Abelian finite
homomorphic image.
%
%
%
Thus we again arrive at a contradiction to property~(i).
Therefore, $G$ is not K\"ahler.
\end{pf*}

\bibliographystyle{amsalpha}

\end{document}